\documentclass[12pt,reqno]{article}
\usepackage{enumerate}
\usepackage{fullpage}
\usepackage{amsmath, amssymb, amsthm, url, graphicx,rotating,tkz-graph,relsize}

\newtheoremstyle{plainsl}%
        {\topsep}
        {\topsep}
        {\slshape} % only non-default setting
        {}
        {\normalfont\bfseries}
        {.}
        { }
        {}

\theoremstyle{plainsl}

\newcommand\cref[1]{Corollary~\ref{cor:#1}}

\newcommand\sqr[2]{{\vbox{\hrule height.#2pt
    \hbox{\vrule width.#2pt height#1pt \kern#1pt
        \vrule width.#2pt}\hrule height.#2pt}}}
\renewcommand\qed{%
        \ifmmode\eqno\sqr53
        \else\nolinebreak\ \hfill\sqr53\medbreak\fi}

\numberwithin{equation}{section}
\begin{document}
\thispagestyle{empty}
\setcounter{page}{1}
\title{On an early paper of Maryam Mirzakhani}
\author{William J. Martin  \\
Department of Mathematical Sciences \\
Worcester Polytechnic Institute \\
Worcester, Massachusetts \\
{\tt martin@wpi.edu}} 

\date{Original: September 16, 2017; Revised:  October 18, 2017, January 24, 2026}
\maketitle

\medskip

\begin{abstract}
Maryam Mirzakhani, the first female (and first Iranian) Fields Medalist, passed away on July 14, 2017 at the age of 40. This short note remembers her 1996 article in % this journal
the Bulletin of the Institute of Combinatorics and its Applications 
and her early years as a mathematician.
\end{abstract}

%%%%%%%%%%%%%%%%%%%%%%%%%%%%%%%%%%%%%%%%%%%%%%%%%%%%%%%%%
\section{We all got the news}
\label{Sec:news}

Maryam Mirzakhani was a brilliant mathematician, being recognized for both her Noether-esque 
accomplishments and her tremendous promise with a Fields Medal at the 2014 ICM in Seoul. 
Dr.\ Mirzakhani was the first female to receive the mathematical community's premier award in 
its 80-year history. She was also the first Iranian Fields Medalist, bringing pride and joy to Iran, 
a nation with a long history of great mathematical advances. 

Born in Tehran on May 3, 1977, Professor Mirzakhani died of metastatic breast cancer on July 14 of this year
at the young age of 40. 
In tribute, Iranian President Hassan Rouhani wrote \cite{chang}, ``The unparalleled excellence of the creative scientist and humble person that echoed Iran's name in scientific circles around the world was a turning point in introducing Iranian women and youth on their way to conquer the summits of pride and various international stages.''

Maryam represented Iran at the International Mathematical Olympiad, earning back-to-back gold medals in 1994 and 1995 with a perfect score in this second competition. After completing her bachelors degree at Sharif University of Technology in 1999, she pursued the PhD at Harvard, writing a dissertation
under the supervision of Curtis McMullen in 2004. The title of her thesis was ``Simple geodesics on hyperbolic surfaces and volume of the moduli space of curves".  After a stint at Princeton, she became a professor at Stanford University, a position Mirzakhani held until her untimely death.

Her citation for the Fields Medal celebrates how Mirzakhani 
``has made stunning advances in the theory of Riemann surfaces and their moduli spaces, and led the way to new frontiers in this area.'' The citation continues, observing that  ``her insights have integrated methods from diverse fields, such as algebraic geometry, topology and probability theory.''

I'd argue that there was a good deal of graph-theoretic intuition behind her work as well.

\section{As for so many, it started with graph coloring}

\begin{figure}
\begin{center}
  \begin{tikzpicture}[style=thick,  % Mirzakhani graph
solidvert/.style={   draw,  circle,   fill=black,  inner sep=1.5pt}]
%solidvert/.style={   draw,  circle,   fill=red,  inner sep=3.5pt}]
\foreach \x/\y in {         0/0, 1/-1,1/0,1/1,
2/0,3/0,
4/-1,4/0,4/1,
5/0,6/0,
7/-1,7/0,7/1,
8/0,9/0,
10/-1,10/0,10/1,
11/0}{
\node[solidvert] (c\x\y)  at (1.2*\x,1.2*\y) {};
\node[solidvert] (ne\x\y) at (1.2*\x+0.6,1.2*\y+0.6) {};
\node[solidvert] (se\x\y) at (1.2*\x+0.6,1.2*\y-0.6) {};
\node[solidvert] (sw\x\y) at (1.2*\x-0.6,1.2*\y-0.6) {};
\node[solidvert] (nw\x\y) at (1.2*\x-0.6,1.2*\y+0.6) {};
\draw (sw\x\y) -- (c\x\y) -- (nw\x\y) -- (sw\x\y) -- (se\x\y) -- (ne\x\y) -- (nw\x\y);
\draw (se\x\y) -- (c\x\y)  -- (ne\x\y);
}
\node[solidvert] (infty) at (6.6 , 5) {};
\draw  (4.6,4) -- (infty) -- (5.6,4); \draw  (7.6,4) -- (infty) -- (8.6,4);  \draw  (5.1,4) -- (infty) -- (8.1,4);
\node (dots) at ( 6.6 , 4.6) {$\infty$};
\node (dots) at ( 6.6 , 4.2) {$\cdots$};
\node (inflab) at ( 6.6 , 5.4) {$L^5$};
%%%
\node (v2lab) at ( -0.6 , -0.9) {$L^2$}; 
\node (v1lab) at ( -0.6 , 0.9) {$L^4$}; 
\node (u1lab) at (  0.37 ,0.0) {$L^1$};
\node (v6lab) at ( 0.3 , -1.9) {$L^2$};  
\node (v5lab) at ( 0.3 , -0.9) {$L^5$};  
\node (v4lab) at ( 0.3 , 0.9) {$L^3$};   
\node (v3lab) at ( 0.3 , 1.9) {$L^4$}; 
\node (lc1lab) at (1.57, -1.2 ) {$L^1$}; 
\node (c1lab) at ( 1.57 , 0.0) {$L^1$}; 
\node (rc1lab) at ( 1.57 , 1.2) {$L^1$};
\node (LL1lab) at ( 2.1 , -1.9) {$L^3$};  
\node (LC1lab) at ( 2.1 , -0.9) {$L^4$};  
\node (RC1lab) at ( 2.1 , 0.9) {$L^2$};   
\node (RR1lab) at ( 2.1 , 1.9) {$L^5$};
\node (uc1lab) at ( 2.77 , 0.0) {$L^1$};
%%%
\node (v1lab) at ( -0.6+3.6 , -0.9) {$L^3$}; 
\node (v1lab) at ( -0.6+3.6 , 0.9) {$L^5$}; 
\node (v1lab) at ( 0.37+3.6 , 0.0) {$L^2$};
\node (v1lab) at ( 0.3+3.6 , -1.9) {$L^3$};  
\node (v1lab) at ( 0.3+3.6 , -0.9) {$L^4$};  
\node (v1lab) at ( 0.3+3.6 , 0.9) {$L^1$};   
\node (v1lab) at ( 0.3+3.6 , 1.9) {$L^5$}; 
\node (lc1lab) at ( 1.57+3.6 , -1.2) {$L^2$}; 
\node (c1lab) at ( 1.57+3.6 , 0.0) {$L^2$}; 
\node (rc1lab) at ( 1.57+3.6 , 1.2) {$L^2$};
\node (LL1lab) at ( 2.1+3.6 , -1.9) {$L^1$};  
\node (LC1lab) at ( 2.1+3.6 , -0.9) {$L^5$};  
\node (RC1lab) at ( 2.1+3.6 , 0.9) {$L^3$};   
\node (RR1lab) at ( 2.1+3.6 , 1.9) {$L^4$};
\node (uc1lab) at ( 2.77+3.6 , 0.0) {$L^2$};
%%%
\node (v1lab) at ( -0.6+7.2 , -0.9) {$L^1$}; 
\node (v1lab) at ( -0.6+7.2 , 0.9) {$L^4$}; 
\node (v1lab) at ( 0.37+7.2 , 0.0) {$L^3$};
\node (v1lab) at ( 0.3+7.2 , -1.9) {$L^1$};  
\node (v1lab) at ( 0.3+7.2 , -0.9) {$L^5$};  
\node (v1lab) at ( 0.3+7.2 , 0.9) {$L^2$};   
\node (v1lab) at ( 0.3+7.2 , 1.9) {$L^4$}; 
\node (lc1lab) at ( 1.57+7.2 , -1.2) {$L^3$}; 
\node (c1lab) at ( 1.57+7.2 , 0.0) {$L^3$}; 
\node (rc1lab) at ( 1.57+7.2 , 1.2) {$L^3$};
\node (LL1lab) at ( 2.1+7.2 , -1.9) {$L^2$};  
\node (LC1lab) at ( 2.1+7.2 , -0.9) {$L^4$};
\node (RC1lab) at ( 2.1+7.2 , 0.9) {$L^1$};   
\node (RR1lab) at ( 2.1+7.2 , 1.9) {$L^5$};
\node (uc1lab) at ( 2.77+7.2 , 0.0) {$L^3$};
%%%
\node (v1lab) at ( -0.6+10.8 , -0.9) {$L^2$}; 
\node (v1lab) at ( -0.6+10.8 , 0.9) {$L^5$};
\node (v1lab) at ( 0.37+10.8 , 0.0) {$L^4$};
\node (v1lab) at ( 0.3+10.8 , -1.9) {$L^2$};  
\node (v1lab) at ( 0.3+10.8 , -0.9) {$L^3$};  
\node (v1lab) at ( 0.3+10.8 , 0.9) {$L^1$};   
\node (v1lab) at ( 0.3+10.8 , 1.9) {$L^5$};
\node (lc1lab) at ( 1.57+10.8 , -1.2) {$L^4$}; 
\node (c1lab) at ( 1.57+10.8 , 0.0) {$L^4$}; 
\node (rc1lab) at ( 1.57+10.8 , 1.2) {$L^4$};
\node (LL1lab) at ( 2.1+10.8 , -1.9) {$L^1$};  
\node (LC1lab) at ( 2.1+10.8 , -0.9) {$L^5$};
\node (RC1lab) at ( 2.1+10.8 , 0.9) {$L^2$};   
\node (RR1lab) at ( 2.1+10.8 , 1.9) {$L^3$};
\node (LC1lab) at ( 2.1+11.7 , -0.9) {$L^1$}; 
\node (RC1lab) at ( 2.1+11.7 , 0.9) {$L^3$}; 
\node (c1lab) at ( 1.57+12.0 , 0.0) {$L^4$};
\end{tikzpicture} 
\end{center}
\caption{The Mirzakhani graph $M$ was discovered in 1996 and first appeared in the journal BICA %this journal 
as a proof 
that there exist planar 3-colorable graphs which are not 4-choosable.
\label{FigM}}
%contain the vertices of the component, $\Delta$, of $\Gamma_b$ which contains $a$.}
\end{figure}
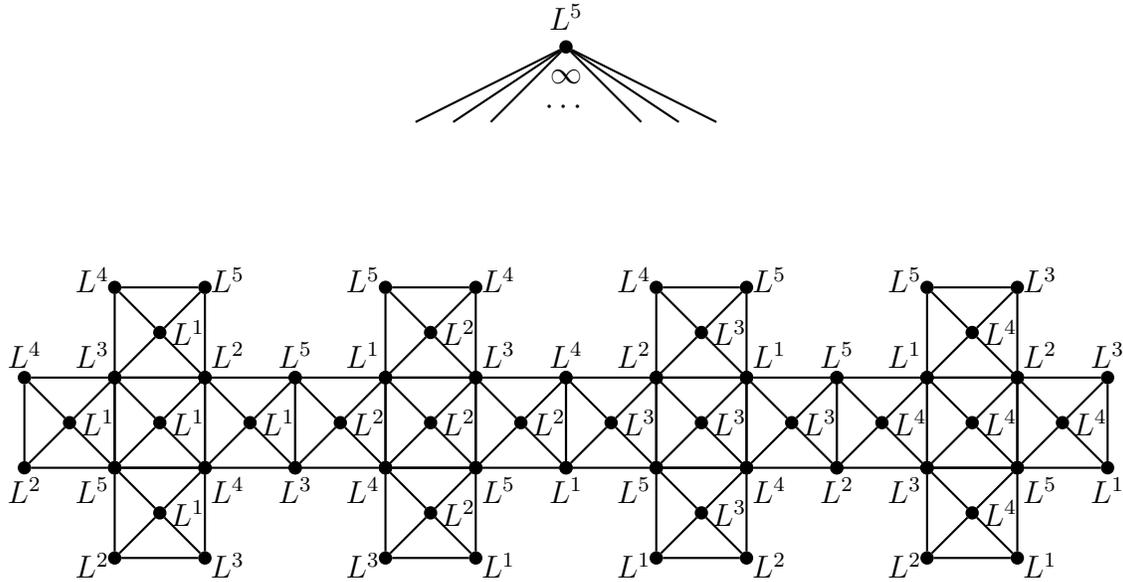

Combinatorics is well-known as a sort of ``gateway drug'' for mathematical research. With problems so easily stated that they can be understood by a casual observer and yet often so challenging that they are impervious
to all but the most keen, persistent and creative puzzle-solvers, graph theory and combinatorics have been
the subject of papers by some of the best mathematicians alive today, even if those mathematicians eventually
settled into other areas of mathematics for their core work.
Of the fifty-six mathematicians having so far been selected for the Fields Medal, twenty-four have published at one
time or another in combinatorics\footnote{MathSciNet primary or secondary classification under 05.}. (The reader may guess which equally enticing mathematical discipline covers another seventeen of the remaining thirty-two.)  

Maryam's paper \cite{mirzakhani} in the BICA\footnote{This was not her first publication. She had already completed one paper \cite{mahmoodian}, also in graph theory, with E.\ S.\ Mahmoodian, which I will discuss below.} was published in 1996, just as she began her undergraduate studies
at Sharif. The title of the paper is ``A small non-4-choosable planar graph''. At the time, list-coloring was a very hot
topic, with Carsten Thomassen's elegant and surprisingly short proof \cite{thomassen} that every planar graph is 5-choosable having just appeared in 1994. Margit Voigt \cite{voigt} had shown in 1993 that there exist planar graphs
which are not 4-choosable. While her proof was quite clever, the example Voigt gave was rather cumbersome: 
the graph has 238 vertices. Maryam's contribution was to reduce the complexity, offering a counterexample 
on 69 vertices which, with a suggestion from her colleague, can be made smaller to give a graph on just 
63 vertices. What is more remarkable about the Mirzakhani example is that it is 3-colorable. So there exist
planar \emph{3-colorable} graphs which are not 4-choosable, disproving a conjecture of Jensen that every planar
3-chromatic graph is 4-choosable.
%. This, therefore, served as a counterexample to another conjecture, by Jensen; 
Alon and Tarsi \cite{alontar} had recently proved that every planar bipartite graph is 3-choosable, so it seemed natural that, for planar graphs at least, the choice number could be bounded above by 
the chromatic number plus one.

Let me briefly review the definitions of these terms. A simple graph $G=(V,E)$ consists of a set $V$ of 
vertices and a set $E$ of edges, each of which joins two distinct vertices. Vertices $u$ and $v$
are adjacent if the edge $uv$ belongs to $E$. A proper coloring of $G$ with colors in set $C$ is a function
$c: V\rightarrow C$ with $c(u) \neq c(v)$ whenever $uv$ is an edge of $G$. Graph $G$ is $k$-colorable
if there is a proper coloring of $G$ into some set $C$ of size $k$. If a subset $L_u$ of $C$ is specified
for each $u\in V$, then a proper list-coloring of $G$ with respect to ``lists''  $\{ L_u \}_{u\in V}$ is
a proper coloring $c: V \rightarrow \cup_u L_u$ satisfying $c(u) \in L_u$ for each vertex $u$. Graph
$G$ is $k$-choosable if such a proper coloring exists for every possible choice of lists $\{ L_u \}_{u\in V}$
such that $|L_u|= k$ for every vertex $u$. The choice number of $G$ is the smallest $k$ for which $G$ is 
$k$-choosable\footnote{Although the graphs we consider here are all finite, it is amusing to note that the Axiom of Choice is the statement that the empty graph (with arbitrary vertex set $V$ and $E=\emptyset$) is list-colorable 
for any assignment of lists provided those lists are all non-empty.}.

%  T. Jensen:  "Such constructions are actually pretty hard to come by and it is a delicate business."

The Mirzakhani $M$ graph is elegant in the way in which its key properties are manifestly evident. We 
reproduce it in Figure \ref{FigM} above, though the arrangement of lists has been slightly modified.
At the top of the diagram is one vertex, which we denote $\infty$, of degree 42. 
Let us call the vertices of degree four in the graph $M-\infty$ ``hubs'' and
those adjacent only to vertices of degree seven in $M-\infty$ ``central hubs''. The graph $M$ is planar and  vertex
$\infty$ is adjacent to every non-hub vertex.

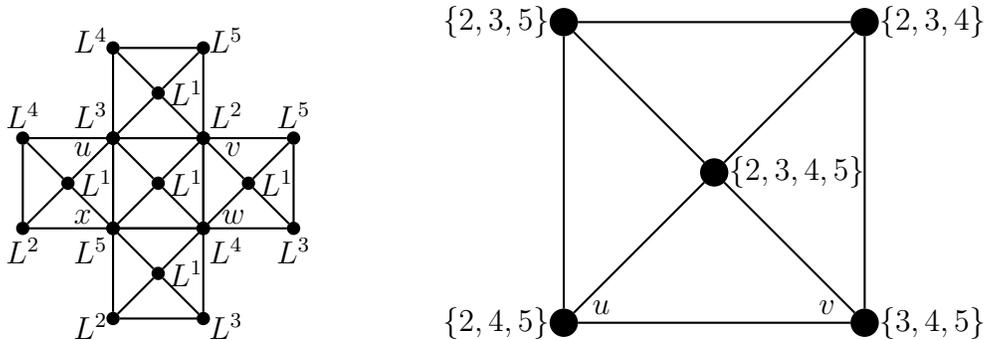
\begin{figure}
\begin{center}
\hfill
  \begin{tikzpicture}[style=thick,  % Mirzakhani graph
solidvert/.style={   draw,  circle,   fill=black,  inner sep=1.5pt}]
%solidvert/.style={   draw,  circle,   fill=red,  inner sep=3.5pt}]
\foreach \x/\y in {         0/0, 1/-1,1/0,1/1,
2/0}{
\node[solidvert] (c\x\y)  at (1.2*\x,1.2*\y) {};
\node[solidvert] (ne\x\y) at (1.2*\x+0.6,1.2*\y+0.6) {};
\node[solidvert] (se\x\y) at (1.2*\x+0.6,1.2*\y-0.6) {};
\node[solidvert] (sw\x\y) at (1.2*\x-0.6,1.2*\y-0.6) {};
\node[solidvert] (nw\x\y) at (1.2*\x-0.6,1.2*\y+0.6) {};
\draw (sw\x\y) -- (c\x\y) -- (nw\x\y) -- (sw\x\y) -- (se\x\y) -- (ne\x\y) -- (nw\x\y);
\draw (se\x\y) -- (c\x\y)  -- (ne\x\y);
}
%%%
\node (v2lab) at ( -0.6 , -0.9) {$L^2$}; 
\node (v1lab) at ( -0.6 , 0.9) {$L^4$}; 
\node (u1lab) at (  0.37 ,0.0) {$L^1$};
\node (v6lab) at ( 0.3 , -1.9) {$L^2$};  
\node (v5lab) at ( 0.3 , -0.9) {$L^5$};  
\node (v4lab) at ( 0.3 , 0.9) {$L^3$};   
\node (v3lab) at ( 0.3 , 1.9) {$L^4$}; 
\node (lc1lab) at (1.57, -1.2 ) {$L^1$}; 
\node (c1lab) at ( 1.57 , 0.0) {$L^1$}; 
\node (rc1lab) at ( 1.57 , 1.2) {$L^1$};
\node (LL1lab) at ( 2.1 , -1.9) {$L^3$};  
\node (LC1lab) at ( 2.1 , -0.9) {$L^4$};  
\node (RC1lab) at ( 2.1 , 0.9) {$L^2$};   
\node (RR1lab) at ( 2.1 , 1.9) {$L^5$};
\node (uc1lab) at ( 2.77 , 0.0) {$L^1$};
%% Name four vertices
\node (u) at  (0.2,0.45) {$u$};
\node (v) at  (2.2,0.45) {$v$};
\node (w) at  (2.2,-0.45) {$w$};
\node (x) at  (0.2,-0.45) {$x$};
%%%
\node (v1lab) at (3 , -0.9) {$L^3$}; 
\node (v1lab) at (3 , 0.9) {$L^5$}; 
\end{tikzpicture} 
%%%%%%%%%%%%%%%%%%%%%%%%%%%%%%%%
\hfill 
%%%%%%%%%%%%%%%%%%%%%%%%%%%%%%%%
  \begin{tikzpicture}[style=thick,  % Mirzakhani graph
solidvert/.style={   draw,  circle,   fill=black,  inner sep=3.5pt}]
%solidvert/.style={   draw,  circle,   fill=red,  inner sep=3.5pt}]
\node[solidvert] (c)  at (0,0) {};
\node[solidvert] (ne) at (2,2) {};
\node[solidvert] (se) at (2,-2) {};
\node[solidvert] (sw) at (-2,-2) {};
\node[solidvert] (nw) at (-2,2) {};
\draw (sw) -- (c) -- (nw) -- (sw) -- (se) -- (ne) -- (nw);
\draw (se) -- (c)  -- (ne);
%%%
\node (clab) at ( 1.1 , 0.0) {$\{2,3,4,5\}$}; % {$L^1$}; 
\node (nelab) at ( 2.9,2) {$\{2,3,4\}$}; 
\node (selab) at ( 2.9,-2) {$\{3,4,5\}$};
\node (swlab) at ( -2.9,-2) {$\{2,4,5\}$};  
\node (nwlab) at ( -2.9,2) {$\{2,3,5\}$};  
\node (u) at  (-1.5,-1.8) {$u$};
\node (v) at  (1.5,-1.8) {$v$};
\end{tikzpicture} 
\hfill $\phantom{X}$
\end{center}
\caption{Subgraphs of the Mirzakhani graph. All lists assigned to the vertices of the first graph 
are 4-element subsets of $\{1,2,3,4,5\}$; e.g.,  $L^1=\{2,3,4,5\}$. Mirzakhani  concisely explained 
why any proper list-coloring of the first graph must use color 1 at some vertex on the outer face. This is seen
by looking at each subgraph isomorphic to the second graph and assuming color 1 is unavailable. 
\label{FigGH}}
%contain the vertices of the component, $\Delta$, of $\Gamma_b$ which contains $a$.}
\end{figure}

The list of colors at each vertex is a subset of $C=\{1,2,3,4,5\}$; the notation $L^j$ denotes the set 
$C \setminus \{j\}$; e.g., $L^4 = \{1,2,3,5\}$. 
Each of the four central hubs, along with the 
four hubs at distance two in $M$ from it has exactly one color from $C$ forbidden; reading left to right, the
lists on the central hubs are $L^1$, $L^2$, $L^3$, $L^4$. If we delete the 20 hub vertices and vertex $\infty$,
we are left with a cycle of length 42; it is therefore evident that $M$ is 3-colorable.  Let us denote this polygon
of length 42 by $P$.

To see that $M$ is not $4$-choosable (cf.~\cite[Exercise 8.4.26]{west}), Mirzakhani first establishes 
properties of the list-colorings of several subgraphs. Consider the first graph in Figure \ref{FigGH}. 
We must convince ourselves that every  proper list-coloring of this graph must use color 1 on the 
outer face. If not, then all lists on the outer face reduce to 
lists of size three. As an example, the subgraph induced on the uppermost hub and its four neighbors 
is reproduced as the second graph in Figure \ref{FigGH}  with its restricted list-coloring problem. Since at most three colors can appear among the four outer vertices, selecting color 5 at the lower left vertex $u$ forces the selection of color 3 at both the upper left and lower right; selecting color 4 at the lower right vertex $v$ likewise forces the selection of color 2 at both the upper right and lower left. 

Allowing for a cyclic permutation $(2453)$ of the colors, we see the same situation at right, bottom and left of the first graph in Figure \ref{FigGH}. We find that $c(u)=5$ forces $c(v)=3$ which, in turn, forces $c(w)=2$ and $c(x)=4$. In fact, any one of these equations forces the three others. Working counterclockwise, the same
is true for the four equations $\{ c(v)=4, c(u)=2, c(x)=3, c(w)=5\}$. The only remaining choice is then
$$c(u)=4, \quad c(v)=5, \quad c(w) = 3, \quad c(x) = 2$$
which also leaves us no choice of color at the central hub. That is, 
in order to afford a color at the central hub, one of the four other hubs must have a neighbor receiving color 1. 

By an appropriate permutation of colors, Mirzakhani is able to string four copies of this graph together 
(each but the last sharing two vertices with the subgraph to its right) in order to obtain the vertex-deleted subgraph $M-\infty$.
By the above argument, color 1 must be used somewhere among the twelve leftmost vertices along polygon $P$. Color 2 must be used on at least one of the twelve vertices consisting of the rightmost two of these and the
ten to their immediate right. Likewise, colors 3 and 4 must be used to color at least one vertex of polygon $P$
in any proper list coloring of $M-\infty$. Since $\infty$ is assigned list $L^5$, this shows that $M$ is not
list-colorable with this assignment of lists. 
%ease with which one may verify its key properties. 

In an email message, Tommy Jensen, co-author of the book ``Graph Coloring Problems'' \cite{jensentoft}
praised the construction, emphasizing the delicacy of finding counterexamples to list-coloring conjectures.
``It is hard
to do this in a short and elegant way, and it seems that
Maryam was able to find a solution that is possibly the most
beautiful of its kind so far'', Jensen wrote.

As a way to remember Dr.\ Mirzakhani with one's students, one might ask those students to study other
fundamental properties of graph $M$. For instance, while $M$ is Hamiltonian and $M-\infty$ has a perfect 
matching, $M-\infty$ is not Hamiltonian.

\section{A community of talented students and their mentors}

Reading Maryam's paper again, one cannot help but be struck by the sense of community in her acknowledgments; she graciously thanks five people, mostly  Olympiad teammates, for their involvement. The problem of finding such a 
counterexample was posed to Maryam by Saieed Akbari who spent eight years training the Iranian Mathematics Olympiad team. Saieed would often pose problems of a combinatorial nature to Maryam and her teammates.
Just as Mirkazhani was entering Sharif University of Technology, Akbari completed his PhD at the same university
under the supervision of  Dr.\ M.\ Mahdavi-Hezavehi and was immediately hired as a faculty member in  the Department of Mathematical Sciences. Saieed had great interest in graph theory (especially algebraic graph 
theory) and combinatorics and was familiar with the construction of Voigt \cite{voigt}. Giving the team only the
necessary definitions and not the history of the problem, Akbari offered 10 US dollars for the construction of
a non-4-choosable planar graph. He notes ``[I]f I had offered 100 US dollars for the construction of such a graph, nobody would dare to attempt the problem.''  After a few days Maryam found him in the corridor and soon 
handed him her one-page proof that there exist planar 3-colorable graphs which are not 4-choosable. In
spite of initial incredulity at the brevity of the solution, Akbari checked the construction over evening tea and 
awarded her the ten-dollar prize the next day. Of her solution, he writes ``I found it to be very creative and beautiful!!''

In my email interview with Saieed Akbari, he recalled Maryam Mirkazhani as polite and disciplined, with a
smile that ``would energize every person that met her.'' He recalls having her in several of his undergraduate and graduate courses. In his linear algebra course, Saieed posed a problem that Maryam solved so elegantly that her 
solution was published in the American Math Monthly \cite{mirmonthly}. He finishes his note to me with the amusing anecdote
that ``Professor Mirzakhani'' once received a request from the US NSF to referee a grant proposal. Apparently,
her impressive achievements masked the fact that she was still an undergraduate student at the time!

\begin{figure}[htb]
\begin{center}
\includegraphics[width=3in]{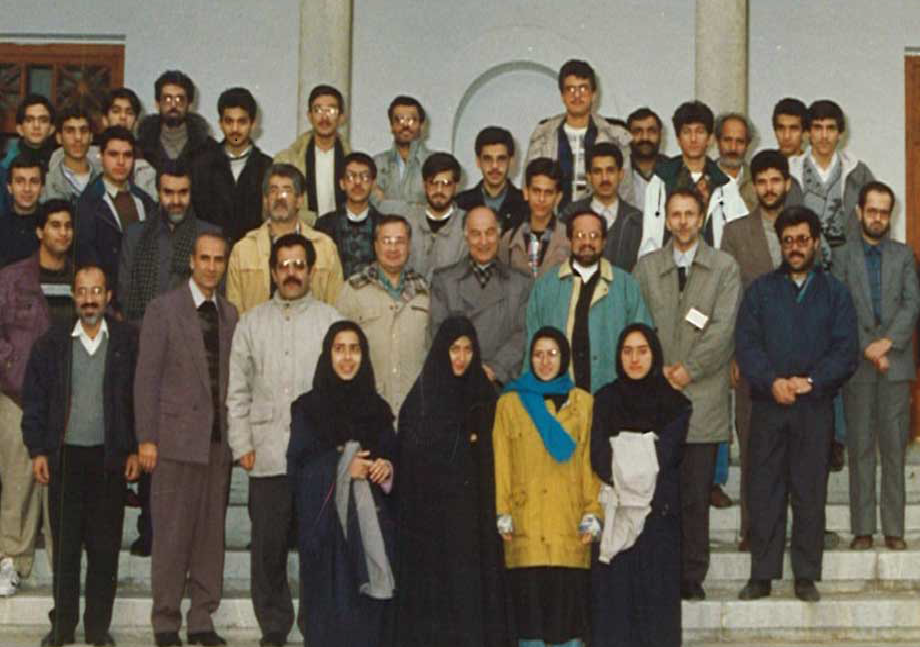}
\end{center}
\caption{Maryam Mirzakhani (front row, right) at the Mathematical Olympiad in Yazd in 1995. Roya Beheshti is also in the front row, at the left. Photo provided by Ebad Mahmoodian. \label{Figthreegirls}}
\end{figure}

As I noted earlier, Maryam's BICA paper was her second publication.  The story of her first publication is also
fascinating.  Ebad Mahmoodian met Maryam when she was 15 years old, having just finished 9th grade. He
and his colleagues at Sharif held summer workshops (or ``math camps'') for students who had completed
their tenth year of high school. They made exceptions for Mirzakhani  and her classmate and friend, Roya Beheshti, who showed such talent in their ninth year. Dr.\ Roya Beheshti Zavareh completed her PhD at MIT in 2003 and is now on the faculty of Washington University in St.\ Louis. In an email interview with Dr.\ Beheshti, I
asked about this experience and Roya responded ``I actually never felt I was at a disadvantage being a female student in the team.'' She explains that it was partly ``because of all the encouragement that we got 
from the professors who were running the math camps''. But Roya added: ``Part of it was because I was with Maryam and it felt empowering to be with another female student, especially someone as strong as her.''
The photograph in 
Figure \ref{Figthreegirls} shows the two of them, both in the front row, at the Iranian Mathematical 
Olympiad in Yazd in 1995.

%Regarding being a female student at the Olympiad team, I actually never felt I was at a disadvantage being a %female student in the team. Part of it was because of all the encouragement that we got 
%from the professors who were running the math camps. Part of it was because I was with Maryam and it felt %empowering to be with another female student specially someone as strong as her.

Having posed the problem of decomposing a 
complete tripartite graph into triangles on the nationwide Olympiad the previous year, Ebad challenged 
the students in
the workshop to find necessary and sufficient conditions for the decomposition of a complete tripartite graph
into 5-cycles. Mahmoodian writes ``Maryam and Roya came up with smart questions. So I invited them to come to the university afterward to continue working on it.''  With the cooperation of her high school principal, Maryam
was able to make visits to the university to collaborate with Ebad, eventually leading to their joint publication
\cite{mahmoodian} which was the first significant contribution to the problem of 5-cycle decompositions of complete tripartite graphs, a problem whose complete solution remains out of reach to this day.

Mahmoodian was so impressed with these young women, that, as a member of the National Mathematics Olympiad Committee, he was able to convince the committee to allow Roya and Maryam to participate 
in their 10th grade. (Ebad notes that, usually, only 11th and 12th graders could participate.) I asked Dr.\ Mahmoodian to share a few words about this amazing woman. Here is part of his response:

\begin{quotation}
{\sl ``I remember Maryam by her  creativity and genuine modesty\ldots.
I am proud of her short fruitful life. Her memory inspires and motivates women all over the world, especially the young female scientists in Iran.''}
\end{quotation}

%%%%%   ********************

Of course, Mirzakhani's main achievements in mathematics were in other areas such as hyperbolic 
geometry \cite{sciam}. She created new tools for the study of moduli spaces of curves and Riemann surfaces; she imported tools from topology and probability theory to Teichm\"{u}ller theory and broke new ground that will remain fruitful for years to come. As an outsider perusing her
papers, I cannot help but think that she still loved the discrete side of things: graphs on surfaces become geodesics, billiard ball trajectories become dynamical systems. Mirzakhani clearly saw mathematics as one.

So many of the best mathematicians start out on problems of a combinatorial nature. But a select few stick with it, or return to that first love. In 2015, already fully into her battle with breast cancer\footnote{As Lamb \cite{sciam} reports, it was already a struggle for Maryam to travel to the ICM in Seoul to receive her Fields Medal. But she kept quiet about the illness.}, %; Jensen describes her as ``very sharp and lively'' that week in South Korea.}, 
Maryam Mirzakhani published, jointly with her spouse and Stanford colleague Jan Vondr\'{a}k, a SODA paper \cite{mv} on Sperner-like colorings of simplicial complexes and their application to fair division.
 This writer can't help but feel that those final months were brutally unfair to Maryam Mirzakhani, her loved ones,  and the community that she inspired so much.

%%%%%%%%%%%%%%%%%%%%%%%%%%%%%%%%%%%%%%%%%%%%%%%%%%%%%%%%%
\section*{Acknowledgments}

Thanks to Saieed Akbari, Roya Beheshti Zavareh, Tommy Jensen and Ebad Mahmoodian for allowing me to interview them via email. Saieed Akbari, Ebad Mahmoodian, Jeff Dinitz and the referee all provided valuable comments and corrections on 
early drafts of the manuscript.

\end{document}